\documentclass[12pt]{amsart}
\usepackage{amsfonts, amssymb,amsmath,amsthm,amscd, graphicx, hyperref}
\usepackage{amsfonts,amssymb,amsmath,amsthm,amscd,graphicx, float}
\newtheorem{theorem}{Theorem}[section]
\newtheorem{lemma}[theorem]{Lemma}
\newtheorem{prop}[theorem]{Proposition}

\theoremstyle{definition}
\newtheorem{definition}[theorem]{Definition}

\theoremstyle{remark}
\newtheorem{remark}[theorem]{Remark}

\numberwithin{equation}{section}

\newcommand{\ot}{\otimes }

\title{Tangle Functors from Semicyclic Representations}

\author{Nathan Druivenga}

\author{Charles Frohman}

\author{Sanjay Kumar}

\address{Department of Mathematics, The University of Iowa, 52242}

\begin{document}

\maketitle

\begin{abstract}  Let $q$ be a $2N$th root of unity where $N$ is odd. Let $U_q(sl_2)$ denote the quantum group with large center corresponding to the lie algebra $sl_2$ with generators $E,F,K$, and $K^{-1}$. A semicyclic representation of $U_q(sl_2)$ is an $N$-dimensional irreducible representation $\rho:U_q(sl_2)\rightarrow M_N(\mathbb{C})$, so that $\rho(E^N)=aId$ with $a\neq 0$, $\rho(F^N)=0$ and $\rho(K^N)=Id$. We construct a tangle functor for framed homogeneous tangles colored with semicyclic representations, and prove that for $(1,1)$-tangles coming from knots, the invariant defined by the tangle functor coincides with Kashaev's invariant. \end{abstract}

\section{Introduction}

Tangle functors have been central in the construction of quantum invariants of links and three-manifolds
since Reshetikhin and Turaev's pioneering work \cite{RT,Tu}. The starting data is a quasitriangular Hopf algebra.  This is a Hopf algebra $H$ along with an invertible element $R$ of a completion of $H\otimes H$ that satisfies three equations. If $R=\sum_i s_i\otimes t_i$ , let $R_{12}=\sum_is_i\otimes t_i \otimes 1$, $R_{13}=\sum_is_i\otimes 1 \otimes t_i$ and $R_{23}=\sum_i 1\otimes s_i\otimes t_i$.  The equations are
 \begin{equation} R \Delta(x)R^{-1}=\Delta'(x),\end{equation} where
$\Delta$ and $\Delta'$ are the comultiplication and the flipped comultiplication in $H$, and
 \begin{equation} (\Delta\otimes 1)R=R_{13}R_{23}, \end{equation} and  \begin{equation}(1\otimes \Delta)(R)=R_{13}R_{12}. \end{equation}
 An element $R\in H\otimes H$ satisfying the above properties is called a {\bf universal $R$-matrix}.

These conditions imply that $R$ satisfies the quantum Yang-Baxter equation,
\begin{equation} R_{12}R_{13}R_{23}=R_{23}R_{13}R_{12}, \end{equation} and
If $\check{R}$ denotes multiplication by $R$ followed by flipping the coordinates, then $\check{R}$ satisfies the standard braid relation,
\begin{equation} \check{R}_{12}\check{R}_{23}\check{R}_{12}=\check{R}_{23}\check{R}_{12}\check{R}_{23}.\end{equation}  Hence quasitriangular Hopf algebras give rise to representations of the braid group.  Tangle functors are the extension of this representation to a representation of the category of framed colored tangles on the category of representations of the quasitriangular Hopf algebra.

The most important examples of quasitriangular Hopf algebras are the Drinfeld-Jimbo deformations of the  universal enveloping algebra of a semisimple Lie algebra  where the deformation parameter is a complex number that is not a root of unity \cite{K}. If there is a universal $R$-matrix, then the representations of the algebra commute under tensor product, that is if $A$ and $B$ are two representations of $H$ then $A\otimes B$ is equivalent to $B\otimes A$.  Roots of unity need to be avoided because the representations of the algebra do not commute under tensor product.  These noncommuting representations are sometimes called {\bf cyclic} representations. The center of the deformed algebra has a subalgebra $Z_0$ of finite index which is a Hopf algebra.  In fact, $Z_0$ is the coordinate ring of the Poisson dual of the Lie group underlying the original Lie algebra. The cyclic representations are associated with points of a branched cover of the dual Lie group \cite{DP,B}.  The cyclic representations are used to construct  {\bf quantum hyperbolic invariants} of three manifolds.  Their restriction to the positive part of the quantum group play a role in Kashaev's invariant\cite{K}, the quantum hyperbolic invariants of Baseilhac and Benedetti \cite{BB}, and implicitly in the representation theory of the Kauffman bracket skein algebra \cite{BL}.

\begin{definition}\label{uq}Let $q\neq 0, \pm 1$ be a complex number. The algebra $U_q(sl_2)$    is generated by $K,K^{-1}$, $E$ and $F$ with relations, 
\begin{equation} \label{eq:rel} KE=q^2EK, \quad KF=q^{-2}FK,  \quad  \mathrm{and} \ EF-FE=\frac{K-K^{-1}}{q-q^{-1}}. \end{equation}
The counit, antipode and comultiplication extend as,
\begin{equation} \label{eq:counit} \epsilon(K)=1, \ \epsilon(E)=0, \ \epsilon(F)=0\end{equation}
\begin{equation} \label{eq:antipode} S(K)=K^{-1}, \ S(E)=-EK^{-1}, \ S(F)=-KF\end{equation}
\begin{equation} \label{comult} \Delta(K)=K\otimes K, \ \Delta(E)=E\otimes K+ 1 \otimes E, \  \Delta(F)=F\otimes 1 + K^{-1}\otimes F.\end{equation} \end{definition}
\smallskip

This is the version of $U_q(sl_2)$ appearing in \cite{DP}, and \cite{B}.
As $U_q(sl_2)$ is not quasitriangular at roots of unity, something must be done in order to get a tangle functor.  In \cite{RT2,KM}, they work with a quotient of $U_q(sl_2)$ where $q=e^{\pi{\bf i}/2r}$ and $r \geq 3$. In this quotient, $K^{2r}=1$ and $E^r=F^r=0$. This algebra is quasitriangular, and the $R$ matrix is,
\begin{equation}\label{ARR}R = q^{\frac{H\otimes H}{2}}\sum_{l=0}^{N-1}\frac{(q-q^{-1})^l}{[l]!}q^{\frac{l(l-1)}{2}}E^l\otimes F^l,\end{equation}  where $K=q^H$. Ohtsuki realized that the requirement that $K^{2r}=1$ was superfluous \cite{Oh}, and the same $R$-matrix defines a tangle functor based on the larger version of the quantum group.

In this paper we find one more extension of the cases where  a tangle functor can be constructed. Here we assume that $q$ is a primitive $2N$th root of unity where $N$ is odd and greater than or equal to $3$. We need $K^N=1$ and $F^N=0$ to define a tangle functor.  We cannot work on the level of the algebra as even when we enforce the relations $K^N=1$ and $F^N=0$ in the algebra, there is no universal $R$-matrix. Hence we work at the level of the representations.   We call representations $\rho$ that satisfy $\rho(K^N)=Id$ and $\rho(F^N)=0$,  {\bf semicyclic}.  First we prove that in these representations a solution to the Quantum Yang-Baxter equation exists and is the image of the $R$-matrix in Equation \ref{ARR}.  Next we prove that for $(1,1)$-tangles the tangle functor returns the same value as if we colored the tangle with the standard irreducible $N$-dimensional representation of $U_q(sl_2)$.  Hitoshi Murakami tells us that Jun Murakami already knew this.  This does not mean that the functor contains no new information. The functor applied   to $(2,2)$-tangles carries  new information  beyond the information in the functor coming from the standard $N$-dimensional representation.

We would like to thank Thomas Kerler who shared with us his insight into a truly elegant proof that the tangle functor based on semicyclic representations returns the same answer as with the standard $N$-dimensional representation.  Given a $(1,1)$-tangle the functor returns a universal polynomial in $K,K^{-1},E,F$ with coefficients in $\mathbb{Z}(q)$, so that the exponential sum of $E$ is equal to the exponential sum of $F$ in every monomial appearing with nonzero coefficient.  Any such monomial can be written as a linear combination of polynomials in $K,K^{-1}$ and the quantum Casimir operator.  The actions of these operators agree in the standard $N$-dimensional representation and the semicyclic representations. Among other things, this means that it is impossible to get quantum hyperbolic invariants from a tangle functor that uses the standard $R$-matrix for $U_q(sl_2)$.

\section{Preliminaries}
Throughout this paper $q=e^{\pi {\bf i}/N}$ where $N\geq 3$ is an odd counting number.
The quantum integer $l$, denoted $[l]$ is defined as,
\[ [l]=\frac{q^l-q^{-l}}{q-q^{-1}}.\]  The quantum factorial is defined recursively by $[0]!=1$, and $[n]!=[n][n-1]!$.  The quantum binomial coefficients are defined by 
\[ \left[\begin{matrix} n \\ k \end{matrix}\right]=\frac{[n]!}{[k]![n-k]!}.\]
The quantum binomial theorem states that if $AB=q^2BA$ then,
\[ (A+B)^n=\sum_{k=0}^nq^{-k(n-k)} \left[\begin{matrix} n \\ k \end{matrix}\right] A^kB^{n-k}.\]

The Weyl algebra $W_q$  is a Hopf algebra generated by $E$, $K$ and $K^{-1}$ with relation $KE=q^2EK$, antipode, $S(K)=K^{-1}$, $S(E)=-EK^{-1}$, counit given by $\epsilon(K)=1$, $\epsilon(E)=0$, and comultiplication $\Delta(K)=K\otimes K$ and $\Delta(E)=E\otimes K+1\otimes E$.  

\subsection{ Cyclic Representations of the Weyl Algebra}\label{Weyl}

Let $V$ be the finite dimensional vector space over the complex numbers with basis $v_i$ where $i\in\{0,1,\ldots, N-1\}$.  Choose a nonzero complex number $a$. Let $M_N(\mathbb{C})$ denote $N\times N$ matrices with complex coefficients, identified with $End(V)$ via the choice of basis.  Define a representation of $W_q$ by,
\[ \rho_{a}: W_q\rightarrow M_N(\mathbb{C}),\]
by $\rho_{a}(K)v_i= q^{2i}v_i$, and $\rho_{a}(E)v_i=v_{i+1}$ when $i<N-1$ and  $\rho_{a}(E)v_{N-1}=av_0$.

Cyclic representations of the Weyl algebra have been studied by Kashaev \cite{K2}, Baseilhac and Benedetti \cite{BB} and Bonahon et al \cite{BL}.

\subsection{ A version of $U_q(sl_2)$}
Let $U_q(sl_2)$ denote the unreduced quantum group as in Definition \ref{uq}. We remind the reader though, that we are only considering the cases where $q=e^{\pi{\bf i}/N}$ where $N\geq 3$ is an odd counting number.

The following equations and proposition will be useful in showing that the $R$-matrix, defined in Section~\ref{sec:Rmatrix}, satisfies the quantum Yang-Baxter equation.
By induction
\[ EF^l=F^lE+[l]F^{l-1}\frac{q^{1-l}K-q^{l-1}K^{-1}}{q-q^{-1}}.\]
Also 
\begin{eqnarray}\label{eq:qH}
q^{\frac{H\otimes H}{2}}(E\otimes 1) q^{-\frac{H\otimes H}{2}} &=&E\otimes K\\
  q^{\frac{H\otimes H}{2}}(1\otimes E) q^{-\frac{H\otimes H}{2}} &=&K\otimes E\\
  q^{\frac{H\otimes H}{2}}(F\otimes 1) q^{-\frac{H\otimes H}{2}} &=&F\otimes K^{-1}\\
 q^{\frac{H\otimes H}{2}}(1\otimes F) q^{-\frac{H\otimes H}{2}} &=&K^{-1}\otimes F
\end{eqnarray}
where $K = q^H$.

\begin{prop}\label{prop:deltaH} $$(\Delta\otimes Id)(q^{\frac{H\otimes H}{2}})=q^{\frac{\Delta(H)\otimes H}{2}} \quad \textrm{and}  \quad
 (Id\otimes \Delta)(q^{\frac{H\otimes H}{2}})=q^{\frac{H\otimes \Delta(H)}{2}}.$$
\end{prop}

 \proof The fact that $\Delta(K)=K\otimes K$ forces 
$\Delta(H) = H \otimes Id + Id \otimes H$.
Using induction, it can be shown that 
$$\Delta(H^n) = \sum_{k=0}^n {{n} \choose {k}} H^k \otimes H^{n-k}.$$

Then,
\begin{eqnarray*}
\left( \Delta \otimes Id \right) \exp{\left( \frac{h}{4} H \otimes H \right) }&=& \left( \Delta \otimes Id \right) \left(\sum_{n=0}^\infty \frac{h^n}{4^n n!} H^n \otimes H^n   \right)\\
&=& \sum_{n=0}^\infty \sum_{k=0}^n  \frac{h^n}{4^n n!}   {{n} \choose {k}} H^k \otimes H^{n-k}   \otimes H^n\\
&=&\sum_{n=0}^\infty \left( \frac{h^n}{4^n n!}  Id \otimes H^n \otimes H^n \right) \left( \sum_{k=0}^n {{n} \choose {k}}    H^k \otimes H^{-k} \otimes Id              \right)\\
&=& \sum_{n=0}^\infty \left( \frac{h^n}{4^n n!}  Id \otimes H^n \otimes H^n \right) \left( H \otimes H^{-1} \otimes Id + Id \otimes Id \otimes Id      \right)^n \\
&=& \sum_{n=0}^\infty  \frac{h^n}{4^n n!}  \left( Id \otimes H \otimes H  + H \otimes Id \otimes H \right)^n \\
&=&\exp{ \left(\frac{h}{4}  \Delta(H) \otimes Id \right) }. \hspace{2in} \qed
\end{eqnarray*}

\subsection{Extending Cyclic Representations of the Weyl Algebra to $U_q(sl_2)$.}\label{sec:extending}

%Assume that $\mu$ is generic so that $\frac{q^{1-2i-k}\mu -q^{-1+2i+k}\mu^{-1}}{q-q^{-1}}$ is not zero unless $q^{1-2i-k}=1$. 

The cyclic representations of $W_q$ defined in Section \ref{Weyl} can be extended to $N$ representations of $U_q(sl_2)$ as follows.  Let $V$ be the  vector space over the complex numbers with basis $v_k$ where $k\in\{0,1,\ldots, N-1\}$. Fix $i \in \{0, 1, ... ,N-1\}$ and $a \in \mathbb{C}-\{0\}$.  Let $\rho_{a,i}(K)v_k= q^{1-N+2(k-i)}v_k$, and $\rho_{a,i}(E)v_k=v_{k+1}$ when $i<N-1$ and  $\rho_{a,i}(E)v_{N-1}=av_0$.  Notice that $\rho_{a,i}(E)^N = a \cdot I_N$ and $\rho_{a,i}(K)^N=I_N$ where $I_N$ is the $N \times N$ identity matrix.  

We find a representation $\rho_{a,i}$ of $F$ such that $\rho_{a,i}(F)^N = 0$.  Under this condition, the equation $\rho_{a,i}(EF-FE)=\rho_{a,i}(\frac{K-K^{-1}}{q-q^{-1}})$ leads to the following solution for $\rho_{a,i}(F)$.    

Denote by $\overline{j}$ the equivalence class of $j \mod N$ which we treat as the remainder on division by $N$. Define

\begin{displaymath}
  \rho_{a,i}(F)v_{\overline{i+k}} = \left\{
     \begin{array}{lr}
       \sum_{j=0}^{\overline{k-1}} -[2(k-j)+N-1]v_{\overline{i+k-1}} & : \overline{i+k} \ne \overline{0} \\
& \\
       \sum_{j=0}^{\overline{k-1}} - \frac{[2(k-j)+N-1]}{a}v_{\overline{i+k-1}} & : \overline{i+k} = \overline{0}
     \end{array}
   \right.
\end{displaymath}

 Equivalently, $\rho_{a,i}(F)$ can be expressed using the sum to product formula,

\begin{displaymath}
  \rho_{a,i}(F)v_{\overline{i+k}} = \left\{
     \begin{array}{lr}
       [k] [N-k] v_{\overline{i+k-1}} & : \overline{i+k} \ne \overline{0} \\
& \\
       \frac{[k] [N-k]}{a}v_{\overline{i+k-1}} & : \overline{i+k} = \overline{0}
     \end{array}
   \right.
\end{displaymath} 
The definition of these $N$ representations leads to the following relations.
\begin{eqnarray}\label{eq:key}
 \rho_{a,i}(F)v_j& =& \rho_{a,i+k}(F)v_{\overline{j+k}}\\
 \rho_{a,i}(K)v_j& =& \rho_{a,i+k}(K)v_{\overline{j+k}}\\
 \rho_{a,i}(E)v_j&= & \rho_{a,i+k}(E)v_j \quad \forall k 
\end{eqnarray}
In fact, all $N$ of the representations are isomorphic via conjugation by $\rho_{a,i}(E)$.  
\begin{eqnarray}
\rho_{a,i}(E)^j \rho_{a,i}(F) \rho_{a,i}(E)^{-j} &=& \rho_{a,\overline{i+j}}(F)\\
\rho_{a,i}(E)^j \rho_{a,i}(K) \rho_{a,i}(E)^{-j} &=& \rho_{a,\overline{i+j}}(K)
\end{eqnarray}

These representations are a subclass of those studied in \cite{GP} and are called {\bf semicyclic representations}. Since there is only one semicyclic representation up to isomorphism, we will work with whichever is most convenient.

\subsection{The Standard Irreducible Representation}

One of the main goals of this paper is to compare the semicyclic and standard  irreducible representations of $U_q(sl_2)$.  Let $\rho_0$ denote the standard $N$-dimensional irreducible representation of  $U_q(sl_2)$.       Fix a basis $\{v_0, v_1, ..., v_{N-1}\}$ of a vector space $V$, then $\rho_0(E), \rho_0(F)$, and $\rho_0(K)$ act on this basis as follows. $ \rho_0(E)v_i = v_{i+1}$ for $i \ne N-1$ and $\rho_0(E)v_{N-1}=0$, $ \rho_0(F)v_i = [i][N-i] v_{i-1}$, $ \rho_0(K)v_i = q^{1-N+2i}v_i $.  The standard representation, $\rho_0$, is similar to the semicyclic representation $\rho_{a,0}$.  In fact, the only difference is $\rho_0(E)v_{N-1}=0$ while $\rho_{a,0}(E)v_{N-1}=av_0$.

\section{ The $R$-matrix}
\label{sec:Rmatrix}

In the case of the standard irreducible representation $\rho_0$, the $R$-matrix  
\[R = q^{\frac{H\otimes H}{2}}\sum_{l=0}^{N-1}\frac{(q-q^{-1})^l}{[l]!}q^{\frac{l(l-1)}{2}}E^l\otimes F^l,\]  satisfies $ R\Delta R^{-1}=\Delta'$ where $\Delta'$ is the flipped comultiplication. 
This equation also holds for  subrepresentations of tensor powers of the two dimensional irreducible representation. This formula first appeared in Kirby and Melvin \cite{KM}, only they were working in the quotient of $U_q(sl_2)$ by $E^N=F^N=0$, and $K^{2N}=1$. This formula for the $R$-matrix also appears in the work of Ohtsuki where he notices that the condition $K^{2N}=1$ is superfluous. It is also used in the
unfolded version of the quantum group, where you add a generator $H$ with $q^{H}=K$.  In the following proposition we will show that the image of $R$ in the representations $\rho_{a,i}$ conjugates the image of comulitplication to the image of flipped comultiplication.

\begin{prop}\label{prop:RdeltaR} For all $i$, $\rho_{a,i}(R)$ and $Z\in U_q(sl_2)$,  satisfies \begin{equation}(\rho_{a,i}\otimes \rho_{a,i})(R\Delta(Z) R^{-1})=(\rho_{a,i}\otimes \rho_{a,i})(\Delta'(Z)). \end{equation}

\end{prop}

\proof In the following proof, we suppress the notation $\rho_{a,i}$ and use $E, F, K$ and $R$ to mean their image under the representation $\rho_{a,i}$. Since $F^N=0$ the proof that $R\Delta(F)-\Delta'(F)R=0$ and $R\Delta(K)-\Delta'(K)R=0$ is the same as the proof appearing in \cite{KM}.  All that is left to show is $R\Delta(E)-\Delta'(E)R=0$.

To simplify notation, write
\[ R=  q^{\frac{H\otimes H}{2}}\sum_{l=0}^{N-1}c_lE^l\otimes F^l.\]
To prove the theorem we show that the equation
\[ R\Delta(E)-\Delta'(E)R=0\] holds for our choice of $R$. 
Using equation~\ref{comult} we get,
\[ q^{\frac{H\otimes H}{2}}(\sum_{l=0}^{N-1}c_lE^l\otimes F^l)(E\otimes K+1\otimes E)-(K\otimes E +E\otimes 1)q^{\frac{H\otimes H}{2}}(\sum_{l=0}^{N-1}c_lE^l\otimes F^l)=0.\]
Using the commutation relations for $q^{\frac{H\otimes H}{2}}$  from equation~\ref{eq:qH}, move it to the front of the second term and factor it out to get,
\[q^{\frac{H\otimes H}{2}}\left((\sum_{l=0}^{N-1}c_lE^l\otimes F^l)(E\otimes K+1\otimes E)-(1\otimes E +E\otimes K^{-1})(\sum_{l=0}^{N-1}c_lE^l\otimes F^l)\right)=0.\]
Cancel $q^{\frac{H\otimes H}{2}}$, distribute and then collect in powers of $E$, to get,
\begin{equation}\label{eqn} \sum_{l=0}^{N-1}c_lE^{l+1}\otimes(F^lK-K^{-1}F^l)+\sum_{l=0}^{N-1}c_lE^l\otimes (F^lE-EF^l)=0.\end{equation}

There is exactly one term with $E^N$ appearing in it,
\[ c_{N-1}E^N\otimes (F^{N-1}K-K^{-1}F^{N-1}).\]
All other terms will cancel as in the proof in \cite{KM} so we only need to show this term is zero.

Under all representations $\rho_{a,i}$,  $F$ has rank $N-1$, so $F^{N-1}$ has rank $1$. Specifically $v_{i-1}$ spans the cokernel of $F^{N-1}$, and its kernel is spanned by all other $v_j$.  Applying $F^{N-1}K-K^{-1}F^{N-1}$ to $v_{i-1}$ under $\rho_{a,i}$  gives
\begin{eqnarray*}
F^{N-1}K-K^{-1}F^{N-1}v_{i-1} &=& F^{N-1}K-q^{2N-2}F^{N-1}K^{-1}v_{i-1}\\
&=& F^{N-1}(K-q^{2N-2}K^{-1})v_{i-1}\\
&=& F^{N-1}( q^{-1-N} - q^{N-1})v_{i-1} \\
&=& 0 \hspace{3in} \qed
\end{eqnarray*}
\begin{remark}Now we can make a standard inductive computation of the $R$-matrix, noting that $c_N=0$.  Going back to Equation \ref{eqn}, commute the $E$ past $F^l$ in the second term, and $K^{-1}$ past $F^{l}$ in the first, and renumber the first term so that $E^{l+1}$ become $E^l$. We get,
\[ \sum_{m=0}^{N-1}c_{m-1}E^m\otimes (F^{m-1}K-q^{2(m-1)}F^{m-1})-\sum_{m=0}^{N-1}c_mE^m\otimes ([m]F^{m-1}\frac{q^{1-m}K-q^{m-1}K^{-1}}{q-q^{-1}})=0.\]
Collecting in powers of $E$ we get the recursive formula,
\[ c_m=\frac{(q-q^{-1})}{[m]}q^{m-1}c_{m-1}.\] Setting $c_0=1$ we get the standard formula for the $R$-matrix given above.  \end{remark}

Recall, a Universal $R$-matrix satisfies  $(\Delta\otimes Id)(R)=R_{13}R_{23}$, and $(Id\otimes \Delta)(R)=R_{13}R_{12}$.  The first of these two equations holds in semicyclic representations because it is true when $E^N=0=F^N$, and $F^N=0$ suffices. This is because the only powers of $E$ and $F$ in the formula that exceed $N-1$ are in fact powers of $F$, for instance see \cite{KM} for the proof.  The second equation does not hold in the semicyclic representations.  To see why, apply Proposition~\ref{prop:deltaH} to the left hand side $(Id\otimes \Delta)(R)$.

\[ (Id \otimes \Delta)(q^{\frac{H\otimes H}{2}}\sum_{l=0}^{N-1}c_lE^l\otimes F^l)=q^{\frac{H\otimes \Delta(H)}{2}}\sum_{l=0}^{N-1}E^l\otimes \Delta(F)^l.\]
Notice that the highest power of $E$ or $F$ that appears is $N-1$.
Since the equation is true when $E^N=F^N=0$ in the algebra it means that
the coefficients of $R_{13}R_{12}$ all monomials where the powers of $E$ and $F$ are less than or equal to $N-1$ agree with the answer above.  For the right hand side, we have
\[ R_{13}R_{12}=q^{\frac{H\otimes 1 \otimes H}{2}}(\sum_{m=0}^{N-1}c_mE^m\otimes 1 \otimes F^m)q^{\frac{H\otimes H \otimes 1}{2}}(\sum_{n=0}^{N-1}c_nE^n\otimes  F^n\otimes 1)\]
Commuting $q^{\frac{H\otimes H \otimes 1}{2}}$ to the front yields,
\[q^{\frac{H\otimes \Delta(H)}{2}}(\sum_{m=0}^{N-1}c_m E^m\otimes K^{-m}\otimes F^m)(\sum_{n=0}^{N-1}c_nE^n\otimes F^n\otimes 1).\]
Now cancel the exponentiated $H\otimes \Delta(H)$ from both sides. We know that all the terms where $m+n\leq N-1$ cancel with the left hand side so the remainder is,
\[ \sum_{m+n\geq N}c_mc_nE^{m+n}\otimes K^{-m}F^n\otimes F^m.\]

The matrices $E^k$ where $k$ ranges from $N$ to $2N-2$ are linearly independent.
Hence each of the parts of the sum where $m+n=k$ needs to be zero.
Letting $k=N$ we get,
\[ \sum_{m=1}^{N-1}c_mc_{n-m}E^N\otimes K^{-m}F^{N-m}\otimes F^m.\]
Consider $v_0\otimes v_{i+1}\otimes v_{i+1}$. The values of this vector under each term of the sum are linearly independent, so the sum is nonzero and the relation $(Id\otimes \Delta)(R)=R_{13}R_{12}$ does not hold under the semicyclic representations.

However, the $R$-matrix still satisfies the Yang-Baxter equation, as long as we are evaluating in the representations $\rho_{a,i}$. \newline

Let $\sigma:A\otimes B\rightarrow B\otimes A$ be the flip $\sigma(Z\otimes W)=W\otimes Z$. Let $P:A\otimes B \otimes C$ denote its extension $P(Z\otimes W\otimes X)=(W\otimes Z\otimes X)$.  If $R=\sum_is_i\otimes t_i$ then
\[ R_{12}=\sum_is_i\otimes t_i\otimes 1, \ R_{13}=\sum_is_i\otimes 1\otimes t_i, \  R_{23}=\sum_i 1\otimes s_i \otimes t_i.\]

\begin{theorem}  Suppose $\rho_{a,i}:U_q(sl_2)\rightarrow M_N(\mathbb{C})$  is any of the semicyclic representations of $U_q(sl_2)$. 
Then
\[ \rho_{a,i}\otimes \rho_{a,i}\otimes \rho_{a,i}(R_{12}R_{13}R_{23}-R_{23}R_{13}R_{12})=0.\] \end{theorem}

\proof The proof follows \cite{KM}. It really only depends on the fact that if we evaluate in $\rho\otimes \phi:U_q(sl_2)\otimes U_q(sl_2)\rightarrow End(V)\otimes M_N(\mathbb{C})$, where $\rho:U_q(sl_2)\rightarrow End(V)$ is an arbitrary representation, and $\phi$ is a semicyclic representation where the kernel of $F$ contains $v_i$, then
\[ (\rho\otimes \phi(R\Delta-\Delta'R)=0,\] and
\[ (\rho\otimes \phi\otimes \psi(\Delta\otimes Id(R)-R_{13}R_{23})=0,\]
when in addition $\psi$ is semicyclic.

To simplify notation we suppress the representation from the formulas.

\[ R_{12}R_{13}R_{23}=R_{12}(\Delta\otimes 1)(R)=(\Delta'\otimes 1)(R) R_{12}=\]
\[P\circ(\Delta\otimes 1)(R)R_{12}=P\circ(R_{13}R_{23})R_{12}=R_{23}R_{13}R_{12}.\]
\qed

\section{The Tangle Functor}
In this section, we define a tangle functor for $(1,1)$-tangles that have been colored with the extended semicyclic representations defined in Section~\ref{sec:extending}.  Let ${e_0,e_1, ... , e_{N-1}}$ be a basis for an $N$-dimensional vector space $V$ over $\mathbb{C}$.  For $v \in V$ and $ \phi \in V^*$ define the cup, cap, and crossing operators as follows. 

\begin{center}
\begin{tabular}{ l  l }
  \raisebox{-3pt}{\includegraphics{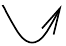}} & $\rightarrow (1\mapsto e_i\otimes e^i)$ \\
  & \\
 \reflectbox{\raisebox{-3pt}{\includegraphics{untwistmin.pdf}}} & $\rightarrow (1\mapsto e^i \otimes K^{-1}(e_i) $ \\
 &\\
  \scalebox{1}[-1]{\raisebox{-8pt}{\includegraphics{untwistmin.pdf}}} & $\rightarrow (\phi \otimes v \mapsto \phi(v)) $ \\
  & \\
  \reflectbox{  \scalebox{1}[-1]{\raisebox{-8pt}{\includegraphics{untwistmin.pdf}}}} & $\rightarrow (v \otimes \phi \mapsto \phi(Kv)) $ \\
  & \\
  
 \reflectbox{ \raisebox{-3pt}{\includegraphics[scale=0.15]{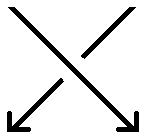}}} & $\rightarrow (e_i \otimes e_j \mapsto \check{R}(e_i\otimes e_j))$ \\   
  & \\
   \raisebox{-3pt}{\includegraphics[scale=0.15]{crossing}} & $\rightarrow (e_i \otimes e_j \mapsto \check{R}^{-1}(e_i\otimes e_j))$ \\
   
\end{tabular}
\end{center}

\begin{theorem}[Turaev \cite{Oh,Tu}]
Two oriented framed tangle diagrams express the same isotopic tangle if and only if the two diagrams are related by a finite sequence of Turaev moves shown below.
\begin{center}
\begin{tabular}{ l  l }
\raisebox{10 pt}{1.} \quad \scalebox{.5}{\raisebox{.1 pt}{\includegraphics{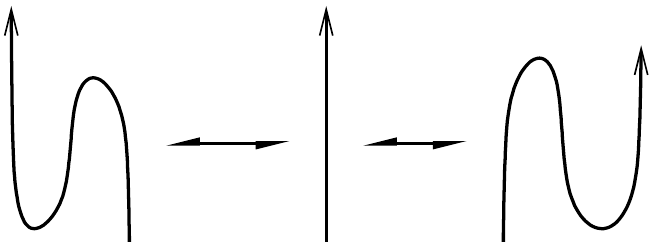}}} & \raisebox{10 pt}{4.}\quad \scalebox{.5}[.37]{\raisebox{8pt}{\includegraphics{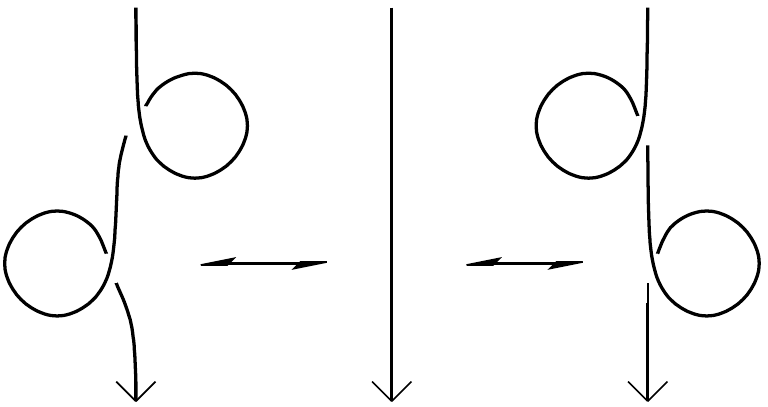}}} \\
  & \\
\raisebox{10 pt}{2.} \quad \scalebox{.5}{\raisebox{\depth}{\scalebox{-1}[-1]{\includegraphics{foist1.pdf}}}} & \raisebox{10 pt}{5.}\quad \scalebox{.5}[.3]{\raisebox{.1pt}{\includegraphics{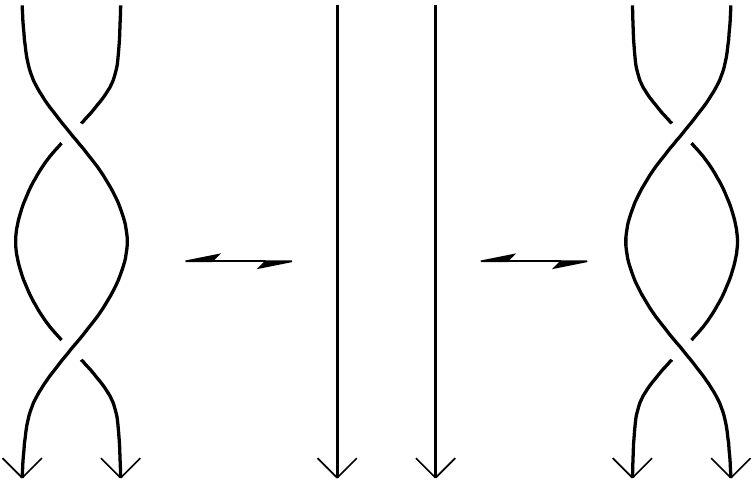}}} \\
 &\\
\raisebox{10 pt}{3.} \quad\scalebox{.45}[.3]{\raisebox{.1pt}{\includegraphics{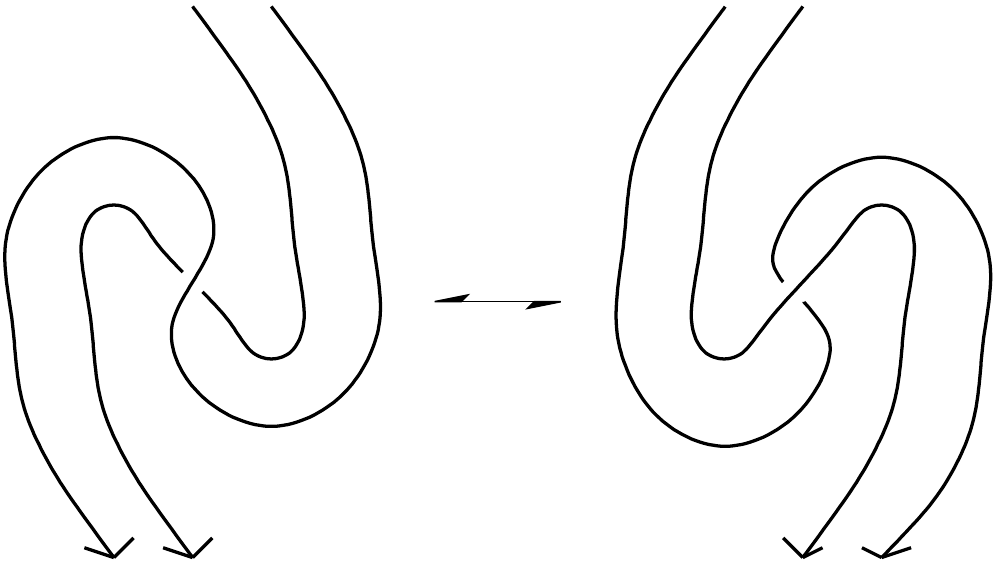}}} &\raisebox{10 pt}{6.} \quad \scalebox{.5}[.3]{\raisebox{.1pt}{   \, \,       \includegraphics{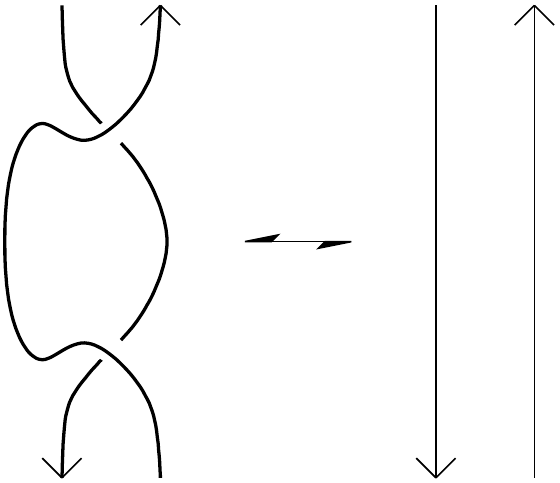}}} \\   
\end{tabular}
\end{center}~\begin{center} \raisebox{10 pt}{7.} \quad \includegraphics[scale=.35]{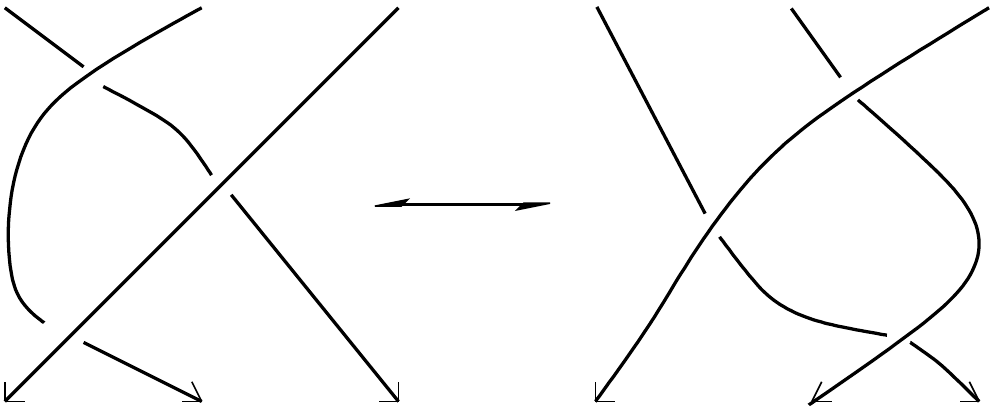} \end{center}
\end{theorem}

We want to show that each of these moves are satisfied using the cap, cup, and crossing definitions given above under specified semicyclic representations.  If so, then we have defined a $(1,1)$ tangle functor. 

\begin{lemma}\label{lem:f_q} Let $q = e^{\frac{\pi i }{N}}$ and $f_q(l) = \frac{(q-q^{-1})^l}{[l]!}q^{\frac{l(l-1)}{2}}$.  The functions $f_q$ and $f_{q^{-1}}$ satisfy,

\begin{equation}f_{q^{-1}}(a-b)f_q(b) = (-1)^{b-a} \frac{(q-q^{-1})^a}{[a-b]![b]!} q^{\frac{a-a^2
}{2}} q^{b(a-1)}. \end{equation}
Also, 
$$\sum_{b=0}^a f_{q^{-1}}(a-b)f_q(b) = 0 =\sum_{b=0}^a f_{q}(a-b)f_{q^{-1}}(b) $$ 
when $0  < a < N$.
\end{lemma}

\proof  The first part is a straightforward algebraic manipulation. The second part is a consequence of the fact that $RR^{-1}=Id\otimes Id$ evaluated in the representation $\rho_{a,i}$. \qed

\begin{prop}\label{prop:tanglefunctor}
The operator $\rho_{a,k}(\check{R})$ for $k \in \{0, 1, ... , N-1\}$ along with the corresponding cup and cap operators arising from $\rho_{a,k}(K)$ define a $(1,1)$-tangle functor. 
\end{prop}
 
 \proof In the following proof, we will assume that $k = \frac{N+1}{2}$ in which case $K(e_i)=q^{2i}e_i$ and $q^{\frac{H \otimes H}{2}}(e_i \otimes e_j) = q^{2ij}(e_i \otimes e_j)$.  Since all the semicyclic representations are isomorphic, we only need to check for one $k$.  Also, we suppress the representation $\rho_{a,k}$ and write $E, F, K$ and $\check{R}$ for their images under $\rho_{a,k}$.  Notice that $1$ and $2$ follow directly from the definition of the cap and cup operators.  Also, $5$ is clearly true since $\check{R} \circ \check{R}^{-1}(e_i \otimes e_j) = e_i \otimes e_j$ and  $\check{R}^{-1} \circ \check{R}(e_i \otimes e_j) = e_i \otimes e_j$.  Since $\check{R}$ satisfies the braid relation, $7$ is satisfied.  It remains to verify 3, 4, and 6.  To simplify notation, we use the modified Einstein notation where there is a sum over each variable that appears at least twice.  
  
 \textbf{Left hand side of 3}:  Here, the arrows are going up so the input to the map is the tensor product of dual vectors $e^i \otimes e^j$.  While passing from the bottom to the top of the diagram, there are two untwisted cups then an application of $\check{R}$ at the crossing, then two untwisted caps.  
 
 \begin{eqnarray*}
e^i \ot e^j &\stackrel{cups}{\mapsto}& e^i \ot e^j \ot  e_k \ot  e_p  \ot e^p \ot e^k     \\
 & \stackrel{\check{R}}{\mapsto}& e^i \ot e^j \ot \check{R}( e_k \ot  e_p)  \ot e^p \ot e^k  \\
 & = &e^i \ot e^j \ot q^{2(p-r)(k+r)} f_q(r) F^r e_p \ot E^r e_k \ot e^p \ot e^k  \\
  & \stackrel{caps}{\mapsto} & q^{2(p-r)(k+r)} f_q(r) (F^r)^j_p (E^r)^i_k e^p \ot e^k \, \, (k = i-r, p = j+r) \\
   & = & q^{2ij}f_q(r) (F^r)^j_{j+r} (E^r)^i_{i-r} e^{j+r} \ot e^{i-r}
\end{eqnarray*}
\textbf{Right hand side of 3}:  As in the case of the left hand side, we calculate the action of the map on $e^i \otimes e^j$.  The only difference is that the cup and cap operators involve the action of $K$.
\begin{eqnarray*}
e^i \ot e^j &\stackrel{cups}{\mapsto}& q^{-2(p+k)} e^p \ot e^k \ot e_k \ot e_p \ot e^i \ot e^j    \\
 & \stackrel{\check{R}}{\mapsto}&  q^{-2(p+k)} e^p \ot e^k \ot \check{R}(e_k \ot e_p) \ot e^i \ot e^j \\
 &=&q^{-2(p+k) + 2(p-r)(k+r)} f_q(r)  e^p \ot e^k \ot F^r e_p \ot E^r e_k \ot e^i \ot e^j \\
 & \stackrel{caps}{\mapsto} & q^{-2(p+k) + 2(p-r)(k+r) + 2(i+j)} f_q(r)  (E^r)^i_k (F^r)^j_p e^p \ot e^k\, \, (k = i-r, p = j+r) \\
   & = & q^{2ij}f_q(r) (F^r)^j_{j+r} (E^r)^i_{i-r} e^{j+r} \ot e^{i-r}
\end{eqnarray*} 
 
\textbf{Left hand side of 4:  }
\begin{eqnarray*}
e_i &\stackrel{cup}{\mapsto}&  q^{-2k} e^k \otimes e_k \otimes e_i \\
     &\stackrel{\check{R}}{\mapsto}&  q^{-2k} e^k \otimes \check{R}(e_k \otimes e_i) \\
     &=&   q^{-2k} q^{-2(i -r)(k +r)} f_q(r) e^k \otimes F^r e_i \otimes E^r e_k \\
     &\stackrel{cap}{\mapsto}& q^{-2k} q^{2(i -r)(k +r)} f_q(r)  (F^r)^k_i  E^r e_k  \quad (k = i-r) \\  
      & \stackrel{cup}{\mapsto} & q^{-2k + 2(i -r)(k +r)}  f_q(r) (F^r)^k_i  E^r e_k \ot e_j \ot e^j \\
      & \stackrel{\check{R}^{-1}}{\mapsto} & q^{-2k + 2(i -r)(k +r)}  f_q(r) (F^r)^k_i \check{R}^{-1} (E^r e_k \ot e_j ) \ot e^j \\
      & = & q^{-2k + 2(i -r)(k +r) -2j(k+r)}  f_q(r)f_{q^{-1}}(s)  (F^r)^k_i E^s e_j \ot F^s E^r e_k \ot e^j \\
      &\stackrel{cap}{\mapsto}& q^{2j -2k + 2(i -r)(k +r) -2j(k+r)}  f_q(r)f_{q^{-1}}(s)  (F^r)^k_i (F^s E^r)^j_k E^s e_j \, \, ( j = i-s) \\
      &=& q^{2(i-s) - 2(i-r) +2i(i-r) - 2i(i-s) }  f_q(r)f_{q^{-1}}(s) (F^{r+s} E^{r+s})^i_i e_i
\end{eqnarray*}

Since $\check{R}$ is an intertwiner, Schur's lemma says the linear map associated to the $(1,1)$ tangle is a multiple of the identity.  Therefore, we can choose any input $0 \le i \le N-1$ to determine the map.  
Set $i=\frac{N+1}{2}$, then $F^r e_i = 0$ for $r > 0$.  This means $r=0$ and $i=k$.  Then $F^s E^r e_k = F^s e_k = 0$ for $s>0$.  Since both $r,s = 0$, the left hand side of 4 is equal to the identity map.  The right hand side of 4 is a similar calculation. \newline
\textbf{Left hand side of 6: }
By the same methods as above, we find that
\begin{eqnarray*}
e^i \ot e_j & \mapsto & \sum_{r,s = 0}^{N-1} q^{-2(r+s)(i-r) +2s} f_{q^{-1}}(r) f_q(s)  (E^{r+s})^i_{i-r-s} (F^{r+s})^{j-r-s}_j e^{i-r-s} \ot e_{j-r-s}
\end{eqnarray*}
Let $a = r+s$ and $b=s$.  Notice that when $a \ge N$, $F^a = 0$ so those terms do not contribute.
\begin{eqnarray*}
e^i \ot e_j & \mapsto & \sum_{b=0}^a \sum_{a=0}^{N-1} q^{-2a(i-a+b) +2b} f_{q^{-1}}(a-b) f_q(b)  (E^{a})^i_{i-a} (F^{a})^{j-a}_j e^{i-a} \ot e_{j-a} \\
 & = & \sum_{a=0}^{N-1} q^{-2a(i-a)} (E^{a})^i_{i-a} (F^{a})^{j-a}_j \left( \sum_{b=0}^a q^{-2b(1-a)} f_{q^{-1}}(a-b) f_q(b) \right)  e^{i-a} \ot e_{j-a}
\end{eqnarray*}
To finish the proof, we need to show that when $a \ne 0$, then $\sum_{b=0}^a q^{-2b(1-a)} f_{q^{-1}}(a-b) f_q(b) = 0$.  Note that when $a = 0$,  $\sum_{b=0}^a q^{-2b(1-a)} f_{q^{-1}}(a-b) f_q(b) = 1$. 
By lemma~\ref{lem:f_q}, when $a \ne 0$, 
\begin{eqnarray*}
\sum_{b=0}^a q^{-2b(1-a)} f_{q^{-1}}(a-b) f_q(b) &=& \sum_{b=0}^a  (-1)^{b-a} \frac{(q-q^{-1})^a}{[a-b]![b]!} q^{\frac{a-a^2}{2}} q^{-b(a-1)} \\
 &=&(-1)^a (q-q^{-1})^a q^{\frac{a-a^2}{2}} \sum_{b=0}^a  (-1)^{b} \frac{ q^{-b(a-1)}   }{[a-b]![b]!}  \\
  & = & (-1)^a q^{a-a^2} \sum_{b=0}^a   f_{q}(a-b) f_{q^{-1}}(b) \\
  & = & 0
\end{eqnarray*}
This concludes the proof of the proposition. \qed

\section{The Invariant of the Figure eight knot}
For notational ease, we continue using the modified Einstein notation where there is a sum over each variable that appears at least twice.  To calculate the invariant, we restrict to the representation $\rho_{a,\frac{N+1}{2}}$.   As in the proof of Proposition~\ref{prop:tanglefunctor} the invariant is calculated by starting at the bottom of the tangle with input vector $e_i$ and applying the proper operator as we pass through cups, caps, and crossings.   Since the number of caps and cups are the same, the final sum will only be over the variables $r,s,t$ and $u$ indicated at each crossing.  Recall that the twisted cup acts as $1 \mapsto e^j \otimes K^{-1}(e_j)$ and the twisted cap as $e^j \otimes e_k \mapsto K(e_k)e^j(e_k)$.  Under $\rho_{a,\frac{N+1}{2}}$, $K^{\pm 1}(e_j) = q^{ \pm 2j}e_j$ and $q^{\frac{H \otimes H}{2}} (e_j \otimes e_k) = q^{2jk}e_j \otimes e_k$.

\begin{figure}[H]
\begin{center}  \includegraphics[scale=0.8]{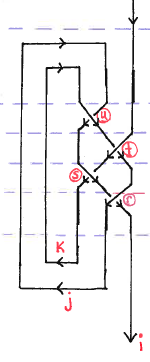}  \caption{ The figure eight knot as a $(1,1)$-tangle}  \end{center}\end{figure}

\hspace{-2.0cm} \begin{tabular}{ r l }
  \small{$e_i $} & \small{$\stackrel{cups}{\mapsto} q^{-2(j+k)} e^j \otimes e^k \otimes e_k  \otimes e_j \otimes e_i $}\\
  & \\ 
   & \small{$\stackrel{\check{R}}{\mapsto}  q^{-2(j+k)}q^{2(i-r)(j+r)} f_q(r) e^j \otimes e^k \otimes e_k  \otimes F^re_i \otimes E^re_j $ } \\
   &\\ 
   & \small{$\stackrel{\check{R}^{-1}}{\mapsto}  q^{-2(j+k)}q^{2(i-r)(j+r)}q^{-2(i-r)k} f_q(r)f_{q^{-1}}(s) e^j \otimes e^k \otimes E^s F^{r}e_i  \otimes F^se_k \otimes E^re_j $}   \\
   &\\
  & \small{$\stackrel{\check{R}}{\mapsto}  q^{-2(j+k)}q^{2(i-r)(j+r)}q^{-2(i-r)k}q^{2(j+r-t)(k-s+t)} f_q(r)f_{q^{-1}}(s)f_q(t) e^j \otimes e^k \otimes E^s F^{r}e_i  \otimes F^t E^re_j \otimes E^t F^{s}e_k $}   \\
   &\\ 
   &\small{$\stackrel{\check{R}^{-1}}{\mapsto}  \overbrace{ q^{-2(j+k)}q^{2(i-r)(j+r)}q^{-2(i-r)k}q^{2(j+r-t)(k-s+t)}q^{-2(i-r+s-u)(j+r-t+u)} f_q(r)f_{q^{-1}}(s)f_q(t)f_{q^{-1}}(u)}^{f(i,j,k,r,s,t,u)}$} \\
   & \qquad \qquad \small{$ e^j \otimes e^k \otimes E^u F^{t} E^re_j  \otimes F^u E^s F^{r}e_i \otimes  E^{t} F^s  e_k $}  \\
   &\\ 
   
   & \small{$\stackrel{caps}{\mapsto} q^{f(i,j,k,r,s,t,u)} \cdot (E^u F^{t} E^r)^k_j (F^uE^s F^{r})^j_i( E^{t} F^s  )^i_k   e_{k-s+t}  \,\, (j=i -r+s-u, k = j +r -t+u)$} \\ 
   & \\
   & \small{${\mapsto} q^{f(i,j(i,r,s,t,u),k(i,r,s,t,u),r,s,t,u)} \cdot (  E^t F^s E^u F^t E^r F^u E^s F^r     )^i_i e_i$}
\end{tabular}
\bigskip 

\textbf{Remark}: In the last step, the exponent of $q$, ${f(i,j,k,r,s,t,u)}$ is replaced with ${f(i,j(i,r,s,t,u),k(i,r,s,t,u),r,s,t,u)}$.  This is due to the linear relations $j = j(i,r,s,t,u)$ and $k=k(i,r,s,t,u)$ coming from the caps.  This will be the case for every knot.  Later we will show that this invariant is the same as Kashaev's invariant so we do not find a closed form for the final solution.

\bigskip
\bigskip

\section{Recovering Kashaev's Invariant}

Given a $(1,1)$-tangle $K$ we denote  the invariant of $K$ coming from the semicyclic representations with $F^N=a$ by $T_{a,N}(K)$ coming from the irreducible semicyclic representations is shown to recover Kashaev's invariant.  Specifically, we show that $T_{a,N}(K)$ evaluated at a $2N^{th}$ root of unity is equivalent to Kashaev's invariant evaluated at a $2N^{th}$ root of unity for any knot.  \newline

A word $w =E^{x_1}F^{x_2} \cdots F^{x_{n-2}} E^{x_{n-1}} F^{x_n}$  is said to be {\bf balanced} if $\sum_{i \ \mathrm{odd}}x_i=\sum_{i \ \mathrm{even}}x_i$. In the prior section, $T_{a,N}(K)$ was calculated for the figure eight knot. That calculation shows that $T_{a,N}(K)$ for any knot will be a sum of balanced  words in $E$ and $F$ with coefficient functions of $q$ to a power.  \newline

Let $E$ be defined as before.  A more general version of $F$ is given by 

   \begin{displaymath}
  (F)^i_j = 
     \begin{cases}
       f_i & : i+1 \equiv j \,  \textrm{mod} \, N\, \textrm{ and } \, i\ne 0 \\
       \frac{f_o}{a}  & :  i = 0 \, \textrm{ and } \, j = N-1\\
       0 & : \textrm{else}
     \end{cases}  
\end{displaymath}               

where $f_i \in \mathbb{C}[q^{\pm1}]$ for $0 \le i \le N-1$. \newline 

Example:  $N = 5$ \newline

$$ F = \begin{pmatrix}  0 & f_1 & 0 & 0 & 0 \\  0 & 0 & f_2 & 0 & 0 \\ 0 & 0 & 0 & f_3 & 0 \\ 0 & 0 & 0 & 0 & f_4 \\ \frac{f_0}{a} & 0 & 0 & 0 & 0   \end{pmatrix}  $$

\begin{prop}\label{prop:noa}
Let $E^{x_1}F^{x_2} \cdots F^{x_{n-2}} E^{x_{n-1}} F^{x_n}$ be balance word, then $ E^{x_1}F^{x_2} \cdots F^{x_{n-2}} E^{x_{n-1}} F^{x_n} $ is a diagonal matrix.  Furthermore, the entries of are elements of $\mathbb{C}[q^{\pm 1}]$.
\end{prop}

\proof  The proof will be by induction on the length of a word $w$ where the total degree of $E$ and $F$ are equal.  For the base case, consider the word $E^m F^m$ of length 2.    For a fixed integer $N$ where $ 0 < m < N $, 

\begin{displaymath}
  (E^m F^m )^i_j = 
     \begin{cases}
       \prod_{k=i}^{\overline{i+m-1}} f_{\overline{N-k}} & : i = j \\
       
       0 & : \textrm{else}
     \end{cases}  
\end{displaymath}    

Also, when $m=0$, $w = id$.  Therefore, for any $0 \le m < N$, the word $w$ satisfies the proposition. \newline

Now assume any word of length less than $n$ satisfies the proposition and let $w = E^{x_1}F^{x_2} \cdots F^{x_{n-2}} E^{x_{n-1}} F^{x_n}$ be a word of length $n$.  There are two cases. \newline

\textbf{Case 1:} $x_n < x_{n-1} $.  In this case, decompose $w$ into the words $a =  E^{x_1}F^{x_2} \cdots F^{x_{n-2}}E^{x_{n-1} -x_n} $  and $b = E^{x_n} F^{x_n}$.  \newline

\textbf{Case 2:}   $x_n > x_{n-1} $.  In this case, it will be necessary to use the following generalized rule to commute the $E^{x_{n-1}}$ past the $F^{x_n}$.

$$  E^c F^d = \sum_{r=0}^c \left( \frac{[c]![d]!}{[r]![c-r]![d-r]!} F^{d-r}E^{c-r} \prod_{k=0}^{r-1} [K;c-d-k]\right) , c < d $$
\newline

where $[K;t] = \frac{Kq^t - K^{-1}q^{-t}}{q - q^{-1}}$.

It is convenient to write the above formula as

$$  E^c F^d = \sum_{r=0}^c  F^{d-r}E^{c-r} D(c,d,r)   , c < d $$

where $D$ is a diagonal matrix comprised entirely of Laurent polynomials in $q$ with no $a$ or $\frac{1}{a}$ appearing.

Applying this formula to the last two terms in the word $w$ gives

\begin{eqnarray*}
w &=& E^{x_1}F^{x_2} \cdots F^{x_{n-2}} E^{x_{n-1}} F^{x_n}\\
  &=& E^{x_1}F^{x_2} \cdots F^{x_{n-2}} \left( \sum_{r=0}^{x_{n-1}}   F^{x_n-r}E^{x_{n-1}-r} D(x_{n-1},x_n,r) \right) \\
  &=&  \sum_{r=0}^{x_{n-1}} E^{x_1}F^{x_2} \cdots F^{x_{n-2}+ x_n - r}   E^{x_{n-1}-r} D(x_{n-1},x_n,r)  \\
\end{eqnarray*}

The result is a word of length less than $n$ which by induction involves no $a's$.  This completes the proof. \qed

The center of $U_q(sl_2)$ is generated by the standard quadratic Casimir defined by $C = EF + \frac{q^{-1} K + q K^{-1}}{(q-q^{-1})^2} = FE +  \frac{qK + q^{-1}K^{-1}}{(q-q^{-1})^2}$.  If the Casimir is evaluated at the standard irreducible representation, $\rho_0$, of $U_q(sl_2)$ then it will clearly not involve any $a's$ since $E^N=0$.  The prior proposition shows that the Casimir will not involve any $a's$ when evaluated at the semicyclic irreducible representations discussed in this paper.

\begin{prop} \label{prop:casimir}
Let $w =E^{x_1}F^{x_2} \cdots F^{x_{n-2}} E^{x_{n-1}} F^{x_n}$ be a balanced word in $E$ and $F$,  then $w$ can be written as a product of the terms $(C - [q^{2\alpha_i}K]_q)$ where $[q^rK]_q = \frac{q^rK-q^{-r}K^{-1}}{(q-q^{-1})^2}$ and $\alpha_i \in \{0,1, ..., N-1\}$.     
\end{prop}

\proof The proof is similar to that of Proposition~\ref{prop:noa} and will be by induction on the length of the word $w$.  Let $w$ be a word of length $2$.  Then $w = E^mF^m$ or $w=F^mE^m$.

By induction and using the formula for the Casimir, it can be shown that  $ E^mF^m = \prod_{i=1}^m(C - [q^{-2(m-i)}K]_{q^{-1}})$ and  $F^mE^m = \prod_{i=1}^m(C - [q^{2(m-i)}K]_{q})$.  With this in hand, the same induction argument from Proposition~\ref{prop:noa} will show that any balanced word can be written in terms of the Casimir and $[q^{2\alpha_i}K]_q$ . \qed

\begin{theorem}
The invariant $T_{a,N}(K)$ evaluated at a $2N$th root of unity is equivalent to Kashaev's invariant, $<K>_N$, evaluated at a $2N$th root of unity for any $(1,1)$-tangle $K$.
\end{theorem}

\proof In this proof, we choose to use $\rho_{a,0}$ when calculating $T_{a,N}(K)$.  Kashaev's invariant is calculated in the same way as $T_{a,N}(K)$ except evaluated at the standard irreducible representation $\rho_0$.  By definition, we have $\rho_{a,0}(K)= \rho_0(K)$.  Therefore, the cup and cap operators are identical under either representation.  Also, $q^{\frac{H \otimes H}{2}}$ acts the same under either representation.  By Proposition~\ref{prop:casimir}, we can complete the proof by checking that $EF$ or $FE$ evaluates the same under both representations.  Notice that $\rho_{a,0}(F)= \rho_0(F)$ and the only difference between $\rho_{a,0}(E)$ and $\rho_0(E)$ is their action on the basis vector $v_{N-1}$.  Thus we only need to check the action of $FE$  on $v_{N-1}$ under each representation. But, under both representations $FEv_{N-1} = 0$.   \qed


\begin{thebibliography}{0000}
\bibitem{B}  Baseilhac, Stephane {\em Quantum coadjoint action and the $6j$-symbols of $U_q(sl2)$.}Interactions between hyperbolic geometry, quantum topology and number theory, 103�143, Contemp. Math., {\bf 541}, Amer. Math. Soc., Providence, RI, 2011. 
\bibitem{BB} 	Baseilhac, Stephane; Benedetti, Riccardo {\em Quantum hyperbolic geometry}, Algebr. Geom. Topol.{\bf  7} (2007), 845–917.
\bibitem{BL}	Bonahon, Francis; Liu, Xiaobo {\em Representations of the quantum Teichmüller space and invariants of surface diffeomorphisms}, Geom. Topol. {\bf 11} (2007), 889–937.
\bibitem{CF} 	Fok, V. V.; Chekhov, L. O. {\em Quantum Teichmüller spaces}, (Russian) Teoret. Mat. Fiz. {\bf 120} (1999), no. 3, 511--528; translation in Theoret. and Math. Phys. {\bf 120} (1999), no. 3, 1245–1259 
\bibitem{DP}  De Concini, C.; Procesi, C. {\em Quantum groups} D-modules, representation theory, and quantum groups (Venice, 1992), 31�140, Lecture Notes in Math., {\bf 1565}, Springer, Berlin, 1993. 
\bibitem{GP} 	Geer, Nathan; Patureau-Mirand, Bertrand {\em G-links invariants, Markov traces and the semicyclic $U_qsl(2)$-modules}, J. Knot Theory Ramifications {\bf 22} (2013), no. 11, 1350063, 28 pp.
\bibitem{K}	Kassel, Christian {\em Quantum groups} Graduate Texts in Mathematics, {\bf 155} Springer-Verlag, New York, 1995. xii+531 pp. ISBN: 0-387-94370-6.
\bibitem{K1}  Kashaev, R. M. {\em Quantization of Teichmüller spaces and the quantum dilogarithm},  Lett. Math. Phys. {\bf 43} (1998), no. 2, 105–115. 
\bibitem{K2}  Kashaev, R. M. {\em Quantum dilogarithm as a 6j-symbol}, Modern Phys. Lett. A {\bf 9} (1994), no. 40, 3757–3768. 
\bibitem{KM} Kirby, Robion; Melvin, Paul {\em The 3-manifold invariants of Witten and Reshetikhin-Turaev for sl(2,C)}, Invent. Math. {\bf 105} (1991), no. 3, 473–545.
\bibitem{Oh} 	Ohtsuki, Tomotada Quantum invariants. A study of knots, 3-manifolds, and their sets. Series on Knots and Everything, 29. World Scientific Publishing Co., Inc., River Edge, NJ, 2002. xiv+489 pp. ISBN: 981-02-4675-7.
\bibitem{R} 	Reshetikhin, N. {\em Quasitriangularity of quantum groups at roots of $1$} Comm. Math. Phys. {\bf 170} (1995), no. 1, 79�99.
\bibitem{RK}  Kashaev, Rinat; Reshetikhin, Nicolai {\em Braiding for quantum $gl_2 $at roots of unity} Noncommutative geometry and representation theory in mathematical physics, 183�197, Contemp. Math., {\bf 391}, Amer. Math. Soc., Providence, RI, 2005. 
\bibitem{RT} 	Reshetikhin, N. Yu.; Turaev, V. G.{\em  Ribbon graphs and their invariants derived from quantum groups} Comm. Math. Phys.{\bf 127} (1990), no. 1, 1–26.
\bibitem{RT2} Reshetikhin, N.; Turaev, V. G. {\em Invariants of 3-manifolds via link polynomials and quantum groups} Invent. Math. {\bf 103} (1991), no. 3, 547–597. 
\bibitem{Tu}	Turaev, V. G. {\em Quantum invariants of knots and 3-manifolds}, de Gruyter Studies in Mathematics, 18. Walter de Gruyter \& Co., Berlin, 1994. x+588 pp. ISBN: 3-11-013704-6 
\end{thebibliography}
\end{document}